\documentclass[12pt]{article}
\usepackage[top=1in,bottom=1in,left=1in,right=1in]{geometry}
\usepackage[english]{babel}
\usepackage{amsthm}
\usepackage{graphicx}
\usepackage{epstopdf}
\usepackage{algorithmic}
\usepackage[version=4]{mhchem}
\usepackage{amsmath}
\usepackage{amssymb}
\usepackage{hyperref}
\usepackage[T1]{fontenc}
\usepackage{soul,xcolor}
\usepackage{cases}
\usepackage{booktabs}
\usepackage{float}
\usepackage{bm}
\usepackage{setspace}
\usepackage{makecell}
\usepackage{longtable}
\usepackage{mathtools}
\usepackage{enumerate}
\usepackage{lineno}
\usepackage{multirow}
\usepackage{pdfpages}
\usepackage{accents}
\usepackage{enumitem}
\usepackage{nomencl}
\makenomenclature
\newcommand*{\dif}{\mathop{}\!\mathrm{d}}

\begin{document}

\title{Decarbonization of Steam Cracking for Clean Olefins Production: Optimal Microgrid Scheduling}

\author{
  Saba Ghasemi Naraghi$^1$ \and Tylee Kareck$^1$ \and Lingyun Xiao$^2$ \and Richard Reed$^3$ \and Paritosh Ramanan$^3$ \and Zheyu Jiang$^{1 \star}$\\
}

\date{
    \normalsize $^1$School of Chemical Engineering, Oklahoma State University, Stillwater, OK 74078\\
    $^2$Operations Research and Industrial Engineering, University of Texas, Austin, TX 78712\\
    $^3$School of Industrial Engineering and Management, Oklahoma State University, Stillwater, OK 74078
}

\maketitle
\vspace{-2em}
\noindent Corresponding author: \texttt{zjiang@okstate.edu}$^\star$

\noindent To appear as a book chapter in Li, C. (Ed.). (2025). \textit{Optimization of Sustainable Process Systems: Multiscale Models and Uncertainties}. Wiley.

\begin{abstract}
    Ethylene is one of the most ubiquitous chemicals and is predominantly produced through steam cracking. However, steam cracking is highly energy- and carbon-intensive, making its decarbonization a priority. Electrifying the steam cracking process is a promising pathway to reduce carbon emissions. However, this is challenged by the intrinsic conflict between the continuous operational nature of ethylene plants and the intermittent nature of renewable energy sources in modern power systems. A viable solution is to pursue a gradual electrification pathway and operate an ethylene plant as a microgrid that adopts diverse energy sources. To optimize the operational strategy of such a microgrid considering uncertainties in renewable energy generation and market prices, in this work, we propose a novel superstructure for electrified steam cracking systems and introduce a stochastic optimization framework for minimizing the operating costs. Results from a case study show that, given the current status of the power grid and renewable energy generation technologies, the process economics and sustainability of electrified steam cracking do not always favor higher decarbonization levels. To overcome this barrier, electricity from the main grid must be cleaner and cheaper, and energy storage costs per unit stored must decrease. Furthermore, it is important for both chemical and power systems stakeholders must seamlessly coordinate with each other to pursue joint optimization in operation.
\end{abstract}

\section{Introduction}

Olefins are hydrocarbons containing one or more double bonds between two adjacent carbon atoms. Olefins are widely used as crucial precursors and essential building blocks in the manufacturing of chemical products, including plastic, detergent, adhesive, rubber, and food packaging. Ethylene is the most important olefin with global annual production exceeding 200 million metric tons \cite{ethyleneproduction}. Like the U.S., the global ethylene market is expected to grow by more than 60\% to 240 billion USD between 2024 and 2034 \cite{ethylenemarket}. Currently, ethylene is almost entirely produced via steam cracking of gaseous and liquid hydrocarbon feedstocks such as ethane, propane, and naphtha. The highly endothermic conversion takes place inside cracking furnaces where preheated feedstock passes through and reacts in furnace coils at very high temperatures (e.g., around $850^{\circ}$C) in a very short residence time (in the order of milliseconds) \cite{ullmann}. The heat is provided by burning the methane fraction byproduct and natural gas in the cracking furnaces. This makes steam cracking one of the most energy and carbon-intensive processes in the chemical industry. Depending on the specific feedstock, it is estimated that 1 to 1.6 tons of \ce{CO2} are released for every ton of ethylene produced \cite{IEA,ullmann}. In fact, steam cracking has been identified by the U.S. Department of Energy as one of the top five energy-intensive refining processes that ``account for the majority of U.S. refining \ce{CO2} emissions and represent the most cost-effective R\&D opportunities to reduce refining emissions'' \cite{roadmap}. Thus, there is a global search for technological advancements and process optimization and intensification to improve the energy efficiency and sustainability of the steam cracking process.

As the U.S. energy landscape continues to transition toward clean, renewable electricity, one promising solution to decarbonize the steam cracking process is to implement electric cracking technology. In 2022, Dow Chemical and Shell started up an experimental electrified cracking unit in Amsterdam, whereas the start-up of a multi-megawatt (MW) pilot plant is expected to take place in 2025 \cite{downews}. Meanwhile, BASF, SABIC, and Linde in partnership have built a demonstration plant in Ludwigshafen, Germany in early 2024 to validate two electric heating concepts of direct (resistive heating directly applied on coils) or indirect (resistive heating elements installed at furnace walls) heating \cite{CEP}. While only 40-45\% of the total firing duty can be transferred through the coils to the process fluids in conventional fired cracking furnaces \cite{purdue}, this number is at least 95\% for electrified cracking furnaces \cite{CEP, maporti2023flexible}. Nevertheless, due to (1) the sheer size of most ethylene plants in the U.S. (whose ethylene production capacities are greater than 1 million tons per year), (2) the need to run these plants around the clock, and (3) the intermittent nature of variable renewable electricity (VRE) from solar and wind whose proportion in the U.S. electricity generation will continue to increase (this number was 14\% in 2023 \cite{eia} and is projected to be 36\% by 2050 \cite{eia2}), it would be economically unrealistic and practically impossible to install massive energy storage systems (e.g., batteries) or perform complete plant reconfiguration to accommodate such a large power demand from electrified crackers \cite{bridge}.

\begin{figure}[ht!]\label{fig_superstructure}
    \centering
    \includegraphics[width = 0.85\textwidth]{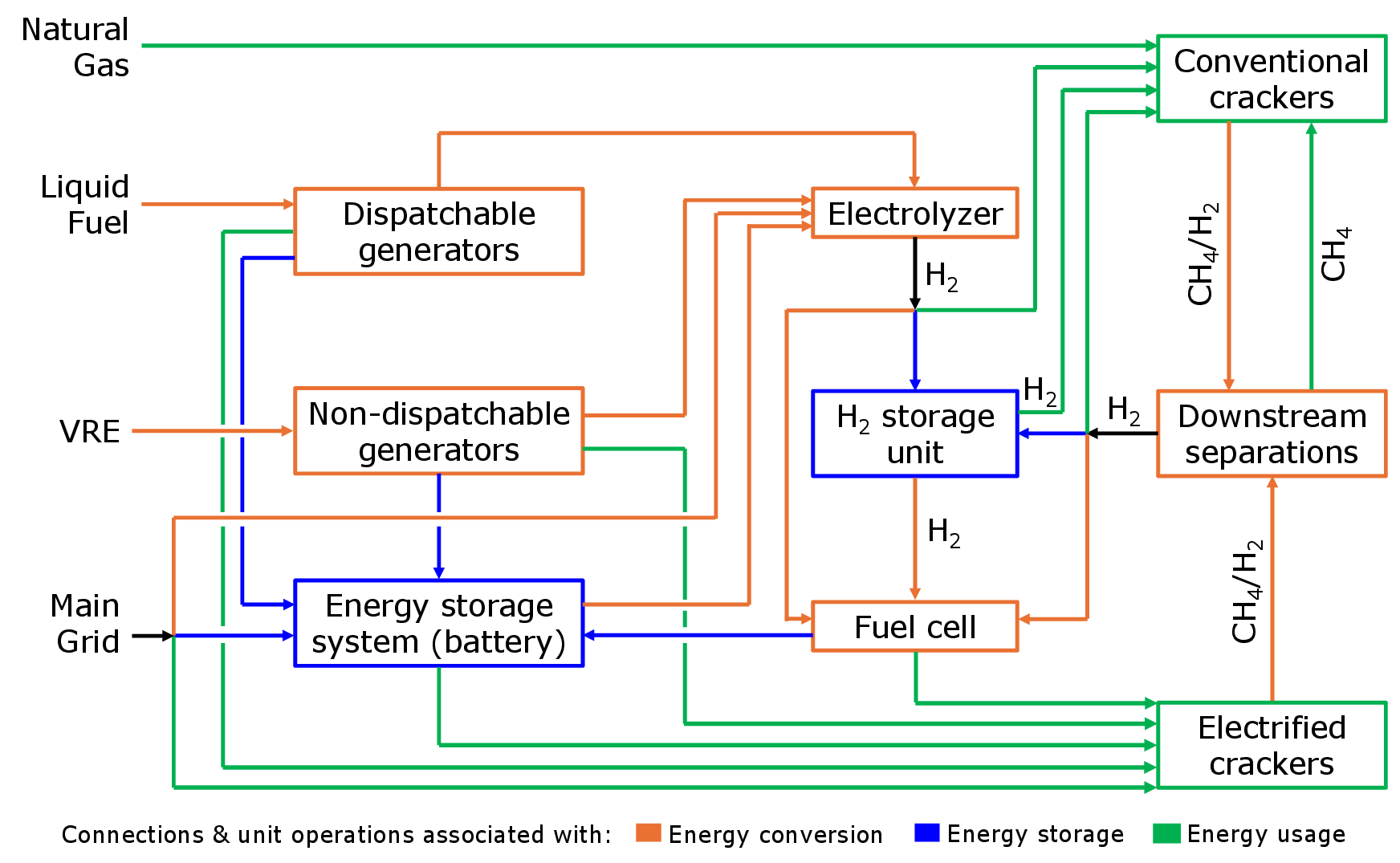}
    \vspace{-1em}
    \caption{Our envisioned framework for using electricity to supply process heat for steam cracking. Diverse energy sources supply heat for both electrified and conventional crackers that are present in the superstructure. Depending on the nature of the energy carriers, the connections shown in the superstructure can represent either energy or mass flows.}
\end{figure}

Accounting for these complications and practical limitations, our vision for using electricity to provide process heat for steam cracking is shown in Figure \ref{fig_superstructure}. We envision that the electrification of steam cracking will take place gradually due to the large capital investment associated with the decommissioning of existing crackers and the installation of new cracker units. Thus, in Figure \ref{fig_superstructure}, both electrified and conventional crackers are present in the superstructure. Battery storage, electrolyzer, and hydrogen storage are used in conjunction with VRE generated onsite to support round-the-clock ethylene plant operation. Electrified crackers can be powered by electricity from the main grid, electricity generated in-house from dispatchable generators and fuel cell units, as well as from batteries. On the other hand, conventional crackers can be powered by fresh natural gas feedstock as well as the methane fraction byproduct (containing \ce{CH4} and \ce{H2}) from both conventional and electrified crackers. Essentially, the future ethylene plant becomes a microgrid, a local electric grid that acts as a single controllable entity with respect to the grid \cite{microgrid}. A microgrid can operate in either grid-connected mode or islanded mode \cite{shahidehpour2012functional,shahidehpour2010role,khodaei2013microgrid,bahramirad2012reliability,flueck2008destination} to account for potential power outages, price fluctuations, and other instabilities associated with the main grid \cite{khodaei2013microgrid,kroposki2008making, kennedy2009reliability,lopez2020reliability,tsikalakis2011centralized,al2021dc}. Therefore, compared with directly plugging into the power system in a centralized manner, microgrid-like steam cracking provides benefits such as improved resilience, economic operation, and flexibility \cite{parhizi2015state}.

This chapter presents the first study of the optimal microgrid scheduling problem for clean olefins production based on the proposed superstructure shown in Figure \ref{fig_superstructure}. We consider a hypothetical ethylene plant located on Texas Gulf Coast. Ethane is chosen to be the major feedstock, which is consistent with most U.S. crackers \cite{CEP}. We consider a deterministic, steady-state operation of the ethylene plant with a fixed total energy demand, ignoring coke formation due to its slow kinetics \cite{ullmann}. On the other hand, the uncertainties associated with VRE generation and market price predictions are considered \cite{gholami2016microgrid}. Note that for this problem, we have a choice of whether to formulation it as a centralized or distributed optimization problem, just like other microgrid scheduling problems studied in the literature \cite{tsikalakis2011centralized,elmouatamid2020review,hatziargyriou2002energy,korpas2006operation,palma2013microgrid}. It turns out that the centralized formulation solves the problem in much less time compared to the distributed formulation based on Benders decomposition. In future works when considering more complex superstructure (e.g., involving downstream processes and plant-wide energy integration pathways), problem specifications (e.g., accounting for coking and decoking), and model formulation (e.g., joint optimization of microgrid and cracker operations), distributed or decentralized formulations may be needed to achieve desired computational performance.

The rest of this chapter is organized as follows. In Section \ref{sec_cracking}, we determine the energy requirement of conventional and electrified cracking by developing a differential-algebraic equation (DAE) numerical model and solving the resulting dynamic optimization problem in \texttt{pyomo.dae} \cite{pyomo}. The energy demand obtained from the mechanistic model is then used to formulate the optimal scheduling problem for the microgrid in Section \ref{sec_formulation}. This scenario-based mixed-integer linear programming (MILP) model is solved in Section \ref{sec_cases}, where several insights are drawn from the case study. Finally, we summarize the results in the concluding remarks in Section \ref{sec_conclusion}.

\section{Dynamic Optimization of Steam Cracking Process}\label{sec_cracking}

The modeling, simulation, and optimization of steam cracking reactors processing different feedstocks have been active research areas for decades following the emergence and growth of petrochemical industries across the world, the recent shale gas boom, and more recently, the industrial decarbonization trend. While process simulator (e.g., Aspen Plus and SPYRO) based modeling and simulation studies are prevalent in the literature \cite{aspen1,aspen2,aspen3,aspen4,aspen5}, mathematical modeling of steam cracking process has also gained considerable research interests \cite{daisonsimulation,fromentbook}, made possible by the pioneering work of Sundaram and Froment \cite{sundaram1977modeling} who developed the first rigorous reaction model. Note that the complete reaction model is based on a free radical mechanism, which is comprehensive and accurate but often results in a stiff nonlinear DAE system due to the large differences in concentration gradient between molecular and free radical species \cite{fromentbook}. Instead, various simplified molecular kinetic models have been proposed \cite{yancheshmeh2013modeling,feli2017investigation,kumarkinetics,onelthesis} to overcome this computational challenge. In this study, we implement an eight-reaction system (see Tables \ref{table_reactions} and \ref{table_kinetics}) previously adopted by \cite{yancheshmeh2013modeling,feli2017investigation} for modeling major ethane cracking reactions.

\begin{table}[ht!]
    \centering
    \begin{tabular}{lll}
        \toprule
        Reaction & Standard enthalpy & Rate expression\\
        & of reaction ($\frac{\text{kJ}}{\text{mol}}$) & \\
        \midrule
        \ce{C2H6 <-> C2H4 + H2} & $ \Delta H^0_1 = 136.33 $ & $r_1 = k_1[\frac{F_{\ce{C2H6}}}{F_{\text{tot}}} (\frac{P}{RT})] - K_{\text{e},1}[\frac{F_{\ce{C2H4}} F_{\text{H}_2}}{F_{\text{tot}}^2} (\frac{P}{RT})^2] $ \\
        \ce{2C2H6 -> C3H8 + CH4} & $ \Delta H^0_2 = -11.56 $ & $r_2 = k_2[\frac{F_{\ce{C2H6}}}{F_{\text{tot}}} (\frac{P}{RT})]$ \\
        \ce{C3H8 -> C3H6 + H2} & $ \Delta H^0_3 = 124.91  $ &$r_3
        = k_3[\frac{F_{\ce{C3H8}}}{F_{\text{tot}}} (\frac{P}{RT})]$ \\
        \ce{C3H8 -> C2H4 + CH4} & $ \Delta H^0_4 = 82.67$ &$r_4
        = k_4[\frac{F_{\ce{C3H8}}}{F_{\text{tot}}} (\frac{P}{RT})]$\\
        \ce{C3H6 <-> C2H2 + CH4} & $ \Delta H^0_5 = 133.45$ &$r_5 = k_5[\frac{F_{\ce{C3H6}}}{F_{\text{tot}}} (\frac{P}{RT})] - K_{\text{e},5}[\frac{F_{\ce{C2H2}} F_{\ce{CH4}}}{F_{\text{tot}}^2} (\frac{P}{RT})^2] $\\
        \ce{C2H2 + C2H4 -> C4H6} & $ \Delta H^0_6 = -171.47$ &$r_6 = k_6[\frac{F_{\ce{C2H2}} F_{\ce{C2H4}}}{F_{\text{tot}}^2} (\frac{P}{RT})^2] $\\ 
        \ce{2C2H6 -> C2H4 + 2CH4} & $ \Delta H^0_7 = 71.10 $ & $r_7 = k_7[\frac{F_{\ce{C2H6}}}{F_{\text{tot}}} (\frac{P}{RT})]$ \\
        \ce{C2H6 + C2H4 -> C3H6 + CH4} & $ \Delta H^0_8 = -22.98 $ & $r_8 = k_8[\frac{F_{\ce{C2H6}} F_{\ce{C2H4}}}{F_{\text{tot}}^2} (\frac{P}{RT})^2] $\\
        \bottomrule
    \end{tabular}
    \caption{Major molecular reactions taking place in ethane cracking.}\label{table_reactions}
\end{table}

\begin{table}[ht!]
    \centering
    \begin{tabular}{lll}
        \toprule
        $k_i = A_i \exp(\frac{-E_i}{RT})$ & Arrhenius constant $A_i$ & Activation energy $E_i$ ($\frac{\mathrm{kJ}}{\mathrm{mol}}$)\\
        \midrule
        $k_1$ & $4.65 \times 10^{13} \, (\frac{\mathrm{m}^3}{\mathrm{mol} \cdot \mathrm{s}})$ & 273.0\\
        $k_2$ & $3.85 \times 10^{11} \, (\frac{1}{\mathrm{s}})$ & 273.0\\
        $k_3$ & $5.89 \times 10^{10} \, (\frac{1}{\mathrm{s}})$ & 215.0\\
        $k_4$ & $4.69 \times 10^{10} \, (\frac{1}{\mathrm{s}})$ & 212.0\\
        $k_5$ & $7.28 \times 10^{12} \, (\frac{\mathrm{m}^3}{\mathrm{mol} \cdot \mathrm{s}})$ & 154.0\\
        $k_6$ & $1.03 \times 10^{9} \, (\frac{1}{\mathrm{s}})$ & 173.0\\
        $k_7$ & $6.37 \times 10^{23} \, (\frac{1}{\mathrm{s}})$ & 530.0\\
        $k_8$ & $7.08 \times 10^{10} \, (\frac{1}{\mathrm{s}})$ & 253.0\\
        $K_{\mathrm{e},1}$ & $8.49 \times 10^{8} \, (\frac{\mathrm{m}^3}{\mathrm{mol} \cdot \mathrm{s}})$ & 136.5\\
        $K_{\mathrm{e},5}$ & $3.81 \times 10^{8} \, (\frac{\mathrm{m}^3}{\mathrm{mol} \cdot \mathrm{s}})$ & 147.2\\
        \bottomrule
    \end{tabular}
    \caption{The rate coefficient expressions in the Arrhenius form.}\label{table_kinetics}
\end{table} 

Following most mathematical models developed in the literature, we model the cracking furnace coils as a one-dimensional plug flow reactor (PFR). This is a reasonable approximation as fluid flow inside the coil is highly turbulent due to the short residence time required to suppress undesired side reactions and coke formation. The configuration specifications of the PFR, which operates in a steady state, are summarized in Table \ref{table_dimensions}. While they are tunable parameters, the default specifications represent common industrial furnaces for ethane cracking. Furthermore, we assume that the key difference between conventional and electrified furnaces lies in the thermal efficiency value $\eta$ \cite{CEP}. Thus, we can estimate the total energy demand of electrified crackers by directly scaling the heat duty results from the dynamic optimization. 

\begin{table}[ht!]
    \centering
    \begin{tabular}{ll}
        \toprule
        Parameter & Default value \\
        \midrule
        Residence time (sec) & 0.2 \\
        Coil length $L$ (m) & 100 \\
        Coil inner diameter $d_\text{t}$ (m) & 0.1 \\
        Inlet steam-ethane mass ratio & 0.3 \\
        $\eta$ of conventional furnace & 0.4 \\
        $\eta$ of electrified furnace & 0.971 \\
        Coil inlet temperature $T_{\text{in}}$ ($^{\circ}$C) & 650 \\
        Coil outlet temperature $T_{\text{out}}$ ($^{\circ}$C) & 850 \\
        Coil inlet pressure $P_{\text{in}}$ (bar) & 3.03 \\
        Coil outlet pressure $P_{\text{out}}$ (bar) & 1.95 \\
        Ethylene yield (kg/kg) & 0.5 \\
        \bottomrule
    \end{tabular}
    \caption{Design and operating specifications of ethane cracking furnace. The ethylene yield, which is the main product specification, is defined as the mass flow ratio between ethane feed and ethylene product.}\label{table_dimensions}
\end{table}

The objective of dynamic optimization is to minimize the total energy demand to achieve an ethylene yield of at least 50\% (which is typical for ethane cracking), subject to mass and energy balances. This is done by adjusting the temperature profile inside the coil, which in turn changes the heat flux profile. While the rigorous momentum balance equation can be incorporated in the model to accurately calculate the pressure drop profile within the coil by accounting for friction losses, P-V-T relationships, and additional losses in the coil bends \cite{fromentbook}, we have identified that incorporating this highly nonlinear ODE does not affect the optimal objective function value, but significantly increases model complexity and computational time. Therefore, we employ a simple analytical expression of Equation \eqref{eqn_pressure} to approximate the pressure drop profile. This equation is derived from the well-known pressure drop expression for flow in pipes \cite{fogler}, namely:
\begin{equation}
    \frac{P(x)}{P_{\text{in}}} = \sqrt{1 - \alpha_p V(x)},
\end{equation}
where $\alpha_p$ is a coil-specific parameter and $V(x) = A_{\text{t}}x$ is the total pipe volume already passed by the fluid.

With this, the dynamic optimization problem is formulated as follows:
\begin{subequations}
    \begin{align}
        \text{minimize \;\;} & \frac{1}{\eta} \int_{0}^{L} q(x) \dif x \\
        \text{subject to:\; } & \frac{\dif F_j}{\dif x} = \frac{\pi d_\text{t}^2}{4} \sum_{i \in I} s_{i,j} r_i, \qquad \forall j \in J \\
        & \sum_{j \in J} F_{j} C_{\text{p},j} \frac{\dif T}{\dif x} = q(x) + \frac{\pi d_{\text{t}}^2}{4} \sum_{i \in I} r_i (- \Delta H_{i}) \\
        & \frac{\dif T }{\dif x} \geq 0 \label{eqn_tempgrad} \\
        & P(x) = P_{\text{in}} \sqrt{1 - \frac{x}{L} \left[1 - \left(\frac{P_{\text{out}}}{P_{\text{in}}} \right) ^2\right]} \label{eqn_pressure}\\
        & \Delta H_{i}(x) =  \Delta H^{0}_{i} + \Delta C_{\text{p},i}(T - 298), \qquad \forall i \in I \\
        & \Delta C_{\text{p},i} = \sum_{j\in J} C_{\text{p},j} s_{i,j}\\
        & F_{\ce{C2H4}}(L) \geq y_{\ce{C2H4}} \frac{\text{MW}_{\ce{C2H6}}}{\text{MW}_{\ce{C2H4}}} F_{\ce{C2H6},\text{in}} \\
        & 0\leq x \leq L;\; T_{\text{in}} \leq T \leq T_{\text{out}}
    \end{align}
\end{subequations}

Here, $I$ and $J$ are sets consisting of all reactions and species included in Table \ref{table_reactions}, respectively. $q(x)$ represents the heat flux (in, e.g., MW), $s_{i,j}$ is the stoichiometry of component $j$ in reaction $i$ ($s_{i,j} > 0$ for products and $<0$ for reactants), molar heat capacity $C_{\text{p},j}$ of pure component $j$ is approximated to be independent of temperature (whose value is obtained at $\frac{T_{\text{in}} + T_{\text{out}}}{2} = 750^\circ$C from Aspen Properties), and $y_{\ce{C2H4}}$ is the ethylene yield specification. Note that by implementing Equation \eqref{eqn_tempgrad}, we ensure that the direction of heat transfer matches with the temperature gradient in the furnace (keep in mind that the flame temperature is as high as 1,150$^\circ$C, which is much higher than the fluid temperature inside the coil). 

This DAE-constrained optimization problem is discretized using an orthogonal collocation method with 13 finite elements and 3 Radau collocation points, resulting in a discretized nonlinear program with 1,863 variables, 1,824 equality constraints, and 39 degrees of freedom. We model the problem in \texttt{pyomo.dae} \cite{pyomo} and find an optimal solution using IPOPT \cite{ipopt} (version 3.14.9) in 34.154 seconds on a Dell Precision 7920 Tower workstation (Intel Xeon Gold 6226R 2.9 GHz, 12$\times$8 GB RAM). The outlet product composition (excluding steam) is shown in Table \ref{table_composition}. The minimum energy requirements for the conventional and electrified crackers are given by 4.27 MWh and 1.75 MWh per ton of ethylene produced. These numbers are used in subsequent sections to study the optimal microgrid scheduling problem.

\begin{table}[ht!]
    \centering
    \begin{tabular}{lll}
        \toprule
        Component & Molar fraction & Mass fraction \\
        \midrule
        Ethane (\ce{C2H6}) & 0.2384 & 0.3753\\
        Ethylene (\ce{C2H4}) & 0.3406 & 0.5\\
        Propane (\ce{C3H8}) & 0.0016 & 0.0038\\
        Propylene (\ce{C3H6}) & 0.0208 & 0.0458\\
        Acetylene (\ce{C2H2}) & 0.0005 & 0.0006\\
        Butadiene (\ce{C4H6}) & 0.0011 & 0.0032\\
        Methane (\ce{CH4}) & 0.0397 & 0.0333\\
        Hydrogen (\ce{H2}) & 0.3573 & 0.0377\\
        \bottomrule
    \end{tabular}
    \caption{Outlet product composition of major species (excluding steam) at the optimal solution obtained by dynamic optimization.} \label{table_composition}
\end{table}

\section{Scenario-Based Optimal Microgrid Scheduling Problem} \label{sec_formulation}

Before formulating the microgrid scheduling problem, we perform a quick calculation \cite{bridge} using the results obtained in the previous section to motivate the need for a diverse portfolio of energy sources. For a plant with 1 million tons/year of ethylene production capacity via electrified cracking, assuming that VRE is on average available for 30\% of a day, then at least 70\% of daily energy (i.e., 3,356 MWh of electricity) needs to be stored for around-the-clock operation of electrified furnaces. Battery storage, based on 100 kWh of the battery pack in a Tesla Model S electric car, would require the battery capacity of 33,560 Tesla Model S cars! Moreover, the actual amount of energy that would need to be stored is likely to be one to two orders of magnitude larger because of daily and seasonal weather variations. This justifies the need to consider a hybrid process heating landscape similar to the superstructure drawn in Figure \ref{fig_superstructure}. Note that in the superstructure, we have both dispatchable and non-dispatchable generators. Dispatchable generators are typically controlled by the microgrid master controllers \cite{khodaei2013microgrid} and are subject to constraints associated with the generation capacity, ramping rates, and minimum on/off time. On the other hand, non-dispatchable generators cannot be controlled by the master controller as their operation depends solely on the availability and capacity of renewable sources. 

In this study, we develop a single-stage scenario-based mixed-integer linear programming (MILP) model to identify the optimal schedule for the microgrid and the optimal degree of electrification for ethane cracking. The sets, parameters, and parameters used in this MILP model are defined below:

\begin{longtable}{ll}
    \textbf{Abbreviations} & \\
    NG & Natural gas\\
    CC, EC & Conventional cracker, electrified cracker\\
    ESS & Energy storage system (battery)\\
    FC & Fuel cell\\
    EL & Electrolyzer\\
    C, DC & Charging, discharging \\
    SU, SD & Startup, shutdown\\
    HS & Hydrogen storage\\
    PV, WT & photovoltaic panel, wind turbine \\
    & \\
    \textbf{Indices} &   \\
    $g$ & Index of dispatchable generator groups $\mathcal{G} = \{1,\dots, G\}$ \\
    $i$ & Index of dispatchable generators belonging to group $g \in \mathcal{G}$, $i \in \mathcal{I}_g$ \\ 
    $t,\, \tau$ & Indices of hourly time periods $\mathcal{T} = \{1,\dots, 24\}$  \\
    $\omega$ & Index of set of scenarios $\Omega$\\
    &  \\
    \textbf{Parameters} &   \\
    $\Delta t$ & Time interval (1 hr) \\
    \(\dot{\text{Q}}_{\text{NG}}\) & Lower heating value of natural gas (13.826 MWh/ton) \\
    \(\dot{\text{Q}}_{\ce{H2}}\) & Lower heating value of hydrogen (33.320 MW/ton) \\
    \(\dot{\text{Q}}_{\ce{CH4}}\) & Lower heating value of methane (13.896 MW/ton) \\
    \(\dot{\text{Q}}_{g}\) & Lower heating value of fuel group $g$ (12.222 MWh/ton) \\
    \(F^{\ce{CH4}/\ce{H2}}_{\text{CC}}\) & Total outlet flow of \ce{CH4} and \ce{H2} from conventional cracker (ton/hr) \\
    \(F^{\ce{CH4}/\ce{H2}}_{\text{EC}}\) & Total outlet flow of \ce{CH4} and \ce{H2} from electrified cracker (ton/hr) \\
    \(f_{\ce{CH4}}\) & Mass fraction of \ce{CH4} in outlet \ce{CH4}/\ce{H2} mixture (0.4692 from Table \ref{table_composition}) \\
    \(r_{\ce{CH4}, \text{sep}}\) & \ce{CH4} recovery in downstream separation units (0.997) \\
    \(r_{\ce{H2}, \text{sep}}\) &  \ce{H2} recovery in downstream separation units (0.99) \\
    \(\text{HSC}\) & Hydrogen storage capacity (ton)\\
    $M^{\ce{H2}}_{\text{HS},\text{start}}$ & Hydrogen storage level at $t=0$\\
    \(\text{ELC} \) & PEM electrolyzer hydrogen production capacity (ton/h) \\
    \(\eta_{\text{EL}}\) & PEM electrolyzer efficiency (0.736) \cite{PEM} \\
    \(\eta_{\text{FC}}\) & Fuel cell efficiency (0.65) \\
    \(\eta_g\) & Generators efficiency (0.60)\\
    \(P^{\ce{H2}}\) & Electrolysis energy requirement with no efficiency loss (39.4 MWh/ton) \\
    \(P_{\text{CC}}\) & Power requirement for conventional cracker (MW) \\
    \(P_{\text{EC}}\) & Power requirement for electrified cracker (MW) \\
    \(\underaccent{\bar}{P}^{\text{FC}},\, \bar{P}^{\text{FC}}\) & Minimum and maximum generation capacity of the fuel cell (MW)\\
    \(\bar{P}^{\text{G}}\) & Maximum power that can be purchased from the grid (MW) \\
    \(P^{\text{WT}}_{t,\omega} \) & Wind power generation (MW) \\
    \(P^{\text{PV}}_{t,\omega} \) & Solar power generation (MW) \\
    \(\underaccent{\bar}{P}^{\text{ESS}}_{\text{C}}, \bar{P}^{\text{ESS}}_{\text{C}} \) & Minimum and maximum ESS (battery) charging power (MW)\\
    \(\underaccent{\bar}{P}^{\text{ESS}}_{\text{DC}}, \bar{P}^{\text{ESS}}_{\text{DC}} \) & Minimum and maximum ESS (battery) discharging power (MW)\\
    \(\text{ESC}\) & Energy storage system capacity (MWh)\\
    $E^\text{ESS}_{\text{start}}$ & Starting ESS charged level at $t=0$\\
    \( \text{MC}^{\text{ESS}},\, \text{MDC}^{\text{ESS}} \) & Minimum battery charging and discharging time (hr) \\
    \(\underaccent{\bar}{P}^{\text{D}}, \bar{P}^{\text{D}}\) & Minimum and maximum generator capacity of dispatchable units (MW)\\
    \(\text{RU}_g^{\text{D}}, \text{RD}_g^{\text{D}}\) & Ramp-up and ramp-down rate for dispatchable group $g$ (MW/hr) \\
    \(\text{UT}_i^{\text{D}}, \text{DT}_i^{\text{D}}\) & Minimum up and down time for dispatchable unit $i$ (h) \\
    \(\text{c}^{\text{Fuel}}_{g} \) & Fuel group $g$ price (\$/ton) \\ 
    \(\text{c}^{\text{G}}_{t,\omega} \) & Electricity market price (\$/MWh) \\
    \(\text{c}^{\text{D}}_{g} \) & Dispatchable units generation cost (\$/MWh)\\
    \(\text{c}^{\text{D}}_{\text{SU},g} \) & Dispatchable units startup cost (\$/MWh)\\
    \(\text{c}^{\text{D}}_{\text{SD},g} \) & Dispatchable units shutdown cost (\$/MWh)\\
    \(\text{c}^{\text{HS}} \) & Hydrogen storage cost (\$/ton)\\
    \(\text{c}^{\text{EL}} \) & PEM electrolyzer operating cost (\$/ton)\\
    \(\text{c}^{\text{FC}} \) & Fuel cell generation cost (\$/MWh)\\
    \(\rho_w\) & Probability of scenario $\omega$ \\
    & \\
    \textbf{Binary variables} &   \\
    \(x^{\text{FC}}_{t, \omega}\) & Fuel cell commitment status at time $t$ and scearnio $\omega$ \\
    \(x^{\text{ESS}}_{\text{C},t,\omega}\) & ESS (battery) charging status at time $t$ and scearnio $\omega$ \\
    \(x^{\text{ESS}}_{\text{DC},t,\omega}\) & ESS (battery) discharging status at time $t$ and scearnio $\omega$ \\
    \(x^{\text{D}}_{i,t,\omega}\) & Dispatchable unit $i$ commitment status at time $t$ and scearnio $\omega$\\
    & \\
    \textbf{Continuous variables} & \textbf{(all under time $t$ and scenario $\omega$)}  \\
    \(F^{\text{NG}}_{\text{CC},t,\omega}\) & Natural gas usage by conventional crackers (ton/h) \\ 
    \(F^{\text{D}}_{i,t,\omega}\) & Fuel usage by generator $i \in \mathcal{I}_g$ (ton/h)\\
    \(F^{\ce{CH4}/\ce{H2}}_{\text{CC},\text{sep},t,\omega}\) & Amount of \(\ce{CH4}\) and \(\ce{H2} \) from CCs to separation unit for energy generation purposes (ton/h)\\
    \(F^{\ce{CH4}/\ce{H2}}_{\text{CC},\text{sep},t,\omega}\) & Amount of \(\ce{CH4}\) and \(\ce{H2} \) from ECs to separation unit for energy generation purposes (ton/h)\\
    \(F^{\ce{CH4}}_{\text{sep},\text{CC},t,\omega}\) & Amount of \ce{CH4} from the separation unit to CCs (ton/h)\\
    \(F^{\ce{H2}}_{\text{sep},\text{CC},t,\omega}\) & Amount of \ce{H2} from the separation unit to CCs (ton/h)\\
    \(F^{\ce{H2}}_{\text{sep},\text{HS},t, \omega}\) & Amount of \ce{H2} from the separation unit to HS (ton/h)\\
    \(F^{\ce{H2}}_{\text{sep},\text{FC},t,\omega}\) & Amount of \ce{H2} from the separation unit to FC (ton/h)\\
    \(F^{\ce{H2}}_{\text{HS},\text{FC},t, \omega}\) & Amount of \ce{H2} from the HS to FC (ton/h)\\
    \(F^{\ce{H2}}_{\text{HS},\text{CC},t, \omega}\) & Amount of \ce{H2} from the HS to CCs (ton/h)\\
    \(M^{\ce{H2}}_{\text{HS},t, \omega}\) & Amount of \ce{H2} stored in HS (ton) \\
    \(F^{\ce{H2}}_{\text{EL},\text{FC},t, \omega}\) & Amount of \ce{H2} produced from the PEM electrolyzer to FC (ton/h)\\
    \(F^{\ce{H2}}_{\text{EL},\text{CC},t, \omega}\) & Amount of \ce{H2} produced from the PEM electrolyzer to CCs (ton/h)\\
     \(F^{\ce{H2}}_{\text{EL},\text{HS},t, \omega}\) & Amount of \ce{H2} produced from the PEM electrolyzer to HS (ton/h)\\
    \(F^{\ce{H2}}_{\text{EL},t, \omega}\) & Amount of \ce{H2} produced from PEM electrolyzer (ton/h) \\
    \(P^{\text{D}}_{\text{EC},g,t,\omega}\) & Amount of power from dispatchable units to ECs (MW)\\
    \(P^\text{G}_{\text{EC},t,\omega}\) & Amount of power from the grid to the ECs (MW)\\
    \(P^\text{ND}_{\text{EC},t,\omega}\) & Amount of power from non-dispatchable units to the ECs (MW)\\
    \(P^\text{ESS}_{\text{EC},t,\omega}\) & Amount of power from the ESS to ECs (MW)\\
    \(P^\text{FC}_{\text{EC},t,\omega}\) & Amount of power from the fuel cell to ECs (MW)\\
     \(P^{\text{D}}_{\text{EL},g,t,\omega} \) & Amount of power from dispatchable units to PEM electrolyzer (MW)\\
     \(P^{\text{G}}_{\text{EL},t,\omega} \) & Amount of power from the grid to PEM electrolyzer (MW)\\
    \(P^{\text{ND}}_{\text{EL},t,\omega} \) & Amount of power from non-dispatchable units to PEM electrolyzer (MW)\\
     \(P^{\text{ESS}}_{\text{EL},t,\omega} \) & Amount of power from ESS to PEM electrolyzer (MW)\\
    \(P^\text{FC}_{\text{ESS},t,\omega}\) & Amount of power from the fuel cell to ESS (MW)\\
    \(P^\text{FC}_{t,\omega}\) & Amount of fuel cell generated power (MW)\\
    \(P^\text{G}_{\text{ESS},t,\omega}\) & Amount of power from the grid to ESS (MW)\\
    \(P^\text{G}_{t,\omega}\) & Amount of power from the grid (MW)\\
    \(P^{\text{ND}}_{\text{ESS},t,\omega} \) & Amount of power from non-dispatchable units to ESS (MW)\\
    \(P^{\text{D}}_{\text{ESS},g,t,\omega}\) & Amount of power from dispatchable units to ESS (MW)\\
    \(P^{\text{ESS}}_{\text{C},t,\omega}\) & Amount of power used to charge the ESS (MW) \\
    \(P^{\text{ESS}}_{\text{DC},t,\omega}\) & Amount of power discharged from ESS (MW)\\
    \(E^{\text{ESS}}_{t,\omega}\) & Amount of stored power in the ESS (MWh)\\
    \(P^{\text{D}}_{i,t,\omega} \) & Amount of power generated by the dispatchable unit $i$ (MW)\\
    \(P^{\text{D}}_{g,t,\omega} \) & Total amount of power generated by dispatchable group $g$ (MW)\\
\end{longtable}

The objective function \eqref{eqn_obj} minimizes the expected total costs, which consist of fuel costs, electricity costs from the grid, energy generation costs, startup and shutdown costs, and hydrogen storage costs. The startup and shutdown costs in the objective function can be easily linearized:
\begin{equation}\label{eqn_obj}
    \begin{aligned}
        \min \quad & \sum_{\omega \in \Omega} \sum_{t \in \mathcal{T}} \rho_{\omega} \Big[\text{c}^{\text{Fuel}}_{\text{NG}} F^{\text{NG}}_{\text{CC},t,\omega} \Delta t +\sum_{g \in \mathcal{G}}\sum_{i \in \mathcal{I}_g} \Big(\text{c}^{\text{Fuel}}_{g} F^{\text{D}}_{i,t,\omega} \Delta t + \text{c}^{\text{D}}_{g} P^{\text{D}}_{i,t,\omega} \\
        & + \text{c}^{\text{D}}_{\text{SU},g} \max \{0, x^{\text{D}}_{i,t,\omega} - x^{\text{D}}_{i,t-1,\omega}\} + \text{c}^{\text{D}}_{\text{SD},g} \max \{0, x^{\text{D}}_{i,t-1,\omega} - x^{\text{D}}_{i,t,\omega} \} \Big)  \\
        & + \text{c}^{\text{G}}_{t,\omega} P^\text{G}_{t,\omega} + \text{c}^{\text{FC}}_{t} P^\text{FC}_{t,\omega} +  \text{c}^{\text{EL}}_{t} F^{\ce{H2}}_{\text{EL},t, \omega} \Delta t+\text{c}^{\text{HS}}_{t} M^{\ce{H2}}_{\text{HS},t, \omega} \Big].
    \end{aligned}
\end{equation}

In terms of the constraints, Equations \eqref{eqn_CC} says that conventional crackers can be powered by natural gas fuel, \ce{CH4} and \ce{H2} produced from cracker units, as well as \ce{H2} coming from the electrolyzer and \ce{H2} storage unit:
\begin{equation}
    \begin{aligned}
        P_{\text{CC}} = & \dot{\text{Q}}_{\text{NG}} F^{\text{NG}}_{\text{CC},t,\omega} + \dot{\text{Q}}_{\ce{CH4}} F^{\ce{CH4}}_{\text{sep},\text{CC},t,\omega} + \\
        & \dot{\text{Q}}_{\ce{H2}} \left(F^{\ce{H2}}_{\text{sep},\text{CC},t,\omega} + F^{\ce{H2}}_{\text{EL},\text{CC},t,\omega} + F^{\ce{H2}}_{\text{HS},\text{CC},t,\omega}\right)
    \end{aligned}
    \qquad \forall t \in \mathcal{T},\, \omega \in \Omega. \label{eqn_CC}
\end{equation}

Equation \eqref{eqn_EC} says that electric crackers can be powered by the main grid, as well as other local fuel-based generators, fuel cells, wind turbines, photovoltaic panels, energy storage systems:
\begin{equation}\label{eqn_EC}
    P_{\text{EC}} = \sum_{g \in \mathcal{G}} P^{\text{D}}_{\text{EC},g,t,\omega} + P^\text{ND}_{\text{EC},t,\omega} + P^\text{G}_{\text{EC},t,\omega} + P^\text{ESS}_{\text{EC},t,\omega} + P^\text{FC}_{\text{EC},t,\omega} \qquad \forall t \in \mathcal{T},\, \omega \in \Omega.
\end{equation}

After steam cracking, the product stream undergoes a series of cooling, drying, compression, and separation operations downstream to obtain individual components from the cracked gas mixture. These steps are required to produce ethylene, propylene, hydrogen, and other value-added hydrocarbons that meet certain specifications. In particular, we are interested in tracking the flows of \ce{CH4} and \ce{H2}, which can be used to power conventional crackers, produce electricity by the fuel cell, or be stored \cite{byproducthydrogen}, as described in Equation \eqref{eqn_gases}. Note that Equation \eqref{eqn_gases} is expressed as inequalities because not all \ce{CH4} and \ce{H2} will be reused.
\begin{equation}
    \begin{split}
        &  F^{\ce{CH4}/\ce{H2}}_{\text{CC},\text{sep},t,\omega}  \leq  F^{\ce{CH4}/\ce{H2}}_{\text{CC}},\\
        &  F^{\ce{CH4}/\ce{H2}}_{\text{EC},\text{sep},t,\omega} \leq F^{\ce{CH4}/\ce{H2}}_{\text{EC}}, \\
        & F^{\ce{CH4}}_{\text{sep},\text{CC},t,\omega} = r_{\ce{CH4},\text{sep}}f_{\ce{CH4}}\left( F^{\ce{CH4}/\ce{H2}}_{\text{CC},\text{sep},t,\omega} +  F^{\ce{CH4}/\ce{H2}}_{\text{EC},\text{sep},t,\omega}\right)\\
        & F^{\ce{H2}}_{\text{sep},\text{CC},t,\omega} + F^{\ce{H2}}_{\text{sep},\text{HS},t,\omega} + F^{\ce{H2}}_{\text{sep},\text{FC},t,\omega} = r_{\ce{H2},\text{sep}}(1-f_{\ce{CH4}})\left(F^{\ce{CH4}/\ce{H2}}_{\text{CC},\text{sep},t,\omega} +  F^{\ce{CH4}/\ce{H2}}_{\text{EC},\text{sep},t,\omega}\right)
    \end{split}
    \qquad \forall t \in \mathcal{T},\, \omega \in \Omega. \label{eqn_gases}
\end{equation}

Equation \eqref{eqn_h2} describes the \ce{H2} mass balances around hydrogen storage unit and electrolyzer:
\begin{equation}
    \begin{split}
        & M^{\ce{H2}}_{\text{HS},t-1, \omega} + \left(F^{\ce{H2}}_{\text{sep},\text{HS},t, \omega} + F^{\ce{H2}}_{\text{EL},\text{HS},t, \omega} -  F^{\ce{H2}}_{\text{HS},\text{FC},t, \omega} - F^{\ce{H2}}_{\text{HS},\text{CC},t,\omega} \right) \Delta t = M^{\ce{H2}}_{\text{HS},t, \omega}, \\
        & F^{\ce{H2}}_{\text{EL},\text{FC},t, \omega} + F^{\ce{H2}}_{\text{EL},\text{HS},t, \omega} + F^{\ce{H2}}_{\text{EL},\text{CC},t, \omega} = F^{\ce{H2}}_{\text{EL},t, \omega}, \\
        & 0 \leq M^{\ce{H2}}_{\text{HS},t, \omega} \leq \text{HSC},\, M^{\ce{H2}}_{\text{HS},0, \omega} = M^{\ce{H2}}_{\text{HS},\text{start}}
    \end{split}
    \qquad \forall t \in \mathcal{T},\, \omega \in \Omega. \label{eqn_h2}
\end{equation}

Equation \eqref{eqn_pem} models hydrogen production in a PEM electrolyzer using electricity from various sources:
\begin{equation}
    \begin{split}
        & F^{\ce{H2}}_{\text{EL},t, \omega} = \frac{\eta_{\text{EL}}}{P^{\ce{H2}}} \left(\sum_{g \in \mathcal{G}} P^{\text{D}}_{\text{EL},g,t,\omega} + P^{\text{ND}}_{\text{EL},t,\omega} + P^\text{G}_{\text{EL},t,\omega} + P^{\text{ESS}}_{\text{EL},t,\omega} \right), \\
        & 0 \leq F^{\ce{H2}}_{\text{EL},t, \omega} \leq \text{ELC},
    \end{split}
    \qquad \forall t \in \mathcal{T},\, \omega \in \Omega. \label{eqn_pem}
\end{equation}

Equation \eqref{eqn_fuelcell} describes constraints related to the fuel cell:
\begin{equation}
    \begin{split}
        & P^{\text{FC}}_{t,\omega} = \eta_{\text{FC}} \dot{\text{Q}}_{\ce{H2}} \left( F^{\ce{H2}}_{\text{HS},\text{FC},t, \omega} + F^{\ce{H2}}_{\text{EL},\text{FC},t, \omega} + F^{\ce{H2}}_{\text{sep},\text{FC},t,\omega} \right), \\
        & \underaccent{\bar}{P}^{\text{FC}} x^{\text{FC}}_{t,\omega} \leq P^{\text{FC}}_{t,\omega} \leq \bar{P}^{\text{FC}} x^{\text{FC}}_{t,\omega},\\
        & P^{\text{FC}}_{t,\omega} = P^{\text{FC}}_{\text{EC},t,\omega} + P^{\text{FC}}_{\text{ESS},t,\omega},
    \end{split}
    \qquad \forall t \in \mathcal{T},\, \omega \in \Omega. \label{eqn_fuelcell}
\end{equation}

The power distribution constraints associated with the main grid are given in Equation \eqref{eqn_grid}:
\begin{equation}
    \begin{split}
        & P^{\text{G}}_{t,\omega} = P^\text{G}_{\text{ESS},t,\omega} + P^\text{G}_{\text{EC},t,\omega} + P^\text{G}_{\text{EL},t,\omega}, \\
        & 0 \leq P^{\text{G}}_{t,\omega} \leq \bar{P}^{\text{G}},
    \end{split}
    \qquad \forall t \in \mathcal{T},\, \omega \in \Omega. \label{eqn_grid}
\end{equation}

Similarly, the power distribution constraints associated with the local non-dispatchable units are shown in Equation \eqref{eqn_nd}:
\begin{equation}
    P^{\text{ND}}_{\text{ESS},t,\omega}  + P^{\text{ND}}_{\text{EC},t,\omega} + P^{\text{ND}}_{\text{EL},t,\omega} = P^{\text{WT}}_{t,\omega} + P^{\text{PV}}_{t,\omega}, \qquad \forall t \in \mathcal{T},\, \omega \in \Omega. \label{eqn_nd}
\end{equation}

For fuel-based local generators, Equation \eqref{eqn_d} considers the power distribution, generation limits, ramping, and up/downtime constraints:
\begin{equation}
    \begin{split}
        & P^{\text{D}}_{i,t,\omega} = \eta_g\dot{\text{Q}}_{g} F^{\text{D}}_{i,t,\omega}, \quad \forall i \in \mathcal{I}_g, \\
        & \underaccent{\bar}{P}^{\text{D}}_{g} x^{\text{D}}_{i,t,\omega} \leq P^{\text{D}}_{i,t,\omega} \leq \bar{P}^{\text{D}}_{g} x^{\text{D}}_{i,t,\omega},\quad \forall i \in \mathcal{I}_g, \\
        & P^{\text{D}}_{g,t,\omega} = \sum_{i \in \mathcal{I}_g} P^{\text{D}}_{i,t,\omega}, \\
        & P^{\text{D}}_{g,t,\omega} = P^{\text{D}}_{\text{ESS},g,t,\omega} + P^{\text{D}}_{\text{EC},g,t,\omega} + P^{\text{D}}_{\text{EL},g,t,\omega}, \\
        & P^{\text{D}}_{i,t,\omega} -P^{\text{D}}_{i,t-1,\omega} \leq \text{RU}^{\text{D}}_g, \quad \forall i \in \mathcal{I}_g, \\
        & P^{\text{D}}_{i,t-1,\omega} -P^{\text{D}}_{i,t,\omega} \leq \text{RD}^{\text{D}}_g, \quad \forall i \in \mathcal{I}_g, \\
        & x^{\text{D}}_{i,\tau,\omega} \geq x^{\text{D}}_{i,t,\omega} - x^{\text{D}}_{i,t-1,\omega}, \quad \forall i \in \mathcal{I}_g,\, \tau \in \{t+1,\dots, \min \{|\mathcal{T}|, t+\text{UT}_{i}-1\}\}, \\
        & 1- x^{\text{D}}_{i,\tau,\omega} \geq x^{\text{D}}_{i,t-1,\omega} - x^{\text{D}}_{i,t,\omega}, \quad \forall i \in \mathcal{I}_g,\, \tau \in \{t+1,\dots, \min \{|\mathcal{T}|, t+\text{DT}_{i}-1\}\},
    \end{split}
    \forall g \in \mathcal{G},\, t \in \mathcal{T},\, \omega \in \Omega. \label{eqn_d}
\end{equation}

The last set of constraints, Equation \eqref{eqn_ess}, considers the power balance and charging/discharging mechanisms of the energy storage system made of batteries. For simplicity, we assume that every battery cell behaves the same.
\begin{equation}
    \begin{split}
        & x^{\text{ESS}}_{\text{C},t,\omega} + x^{\text{ESS}}_{\text{DC},t,\omega} \leq 1, \\
        & x^{\text{ESS}}_{\text{C},\tau,\omega} \geq x^{\text{ESS}}_{\text{C},t,\omega} - x^{\text{ESS}}_{\text{C},t-1,\omega}, \quad \forall \tau \in \{t+1,\cdots, \min \{|\mathcal{T}|, t+\text{MC}^{\text{ESS}}-1\} \}, \\
        & x^{\text{ESS}}_{\text{DC},\tau,\omega} \geq x^{\text{ESS}}_{\text{DC},t,\omega} - x^{\text{ESS}}_{\text{DC},t-1,\omega}, \quad \forall \tau \in \{t+1,\cdots, \min \{|\mathcal{T}|, t+\text{MDC}^{\text{ESS}}-1\} \}, \\
        & \underaccent{\bar}{P}^{\text{ESS}}_{\text{C}} x^{\text{ESS}}_{\text{C},t,\omega} \leq P^{\text{ESS}}_{\text{C},t,\omega} \leq \bar{P}^{\text{ESS}}_{\text{C}} x^{\text{ESS}}_{\text{C},t,\omega}, \\
        & \underaccent{\bar}{P}^{\text{ESS}}_{\text{DC}} x^{\text{ESS}}_{\text{DC},t,\omega} \leq P^{\text{ESS}}_{\text{DC},t,\omega} \leq \bar{P}^{\text{ESS}}_{\text{DC}} x^{\text{ESS}}_{\text{DC},t,\omega}, \\
        & P^{\text{ESS}}_{\text{C},t,\omega} = P^\text{G}_{\text{ESS},t,\omega} + \sum_{g \in \mathcal{G}} P^{\text{D}}_{\text{ESS},g,t,\omega} + P^{\text{ND}}_{\text{ESS},t,\omega} + P^{\text{FC}}_{\text{ESS},t,\omega}, \\
        & P^{\text{ESS}}_{\text{DC},t,\omega} = P^{\text{ESS}}_{\text{EC},t,\omega} + P^{\text{ESS}}_{\text{EL},t,\omega}, \\
        & E^{\text{ESS}}_{t,\omega} = E^{\text{ESS}}_{t-1,\omega} + \left(P^{\text{ESS}}_{\text{C},t,\omega} - P^{\text{ESS}}_{\text{DC},t,\omega}\right) \Delta t, \\
        & 0 \leq E^{\text{ESS}}_{t,\omega} \leq \text{ESC},\, E^{\text{ESS}}_{0,\omega} = E^\text{ESS}_{\text{start}},
    \end{split}
    \quad \forall t \in \mathcal{T},\, \omega \in \Omega. \label{eqn_ess}
\end{equation}

\section{Illustrative Case Studies}\label{sec_cases}

\subsection{Problem Setting}

In this section, the proposed scenario-based MILP model is applied to a hypothetical ethylene plant for illustrative case studies. The plant has an ethylene production capacity of 1 million tons/year. Ethane cracking furnaces in the U.S. have capacities ranging from 100,000 to 250,000 tons of ethylene/year. Thus, we assume the plant has five ethane crackers with 200,000 tons of ethylene/year capacity. By solving this optimal scheduling problem, we hope to provide the first quantitative insights about how steam cracking electrification shall be conducted economically and sustainably. The features of dispatchable generators, energy storage systems, fuel cells, and hydrogen storage are summarized in Table \ref{table_parameters}.

\begin{table}[ht]
    \centering 
    \begin{tabular}{cccccc}
    \toprule
    Fuel group & \# of & Operating cost & Min-max & Min up/down & Ramp-up/down\\
    & units & (\$/MWh) & capacity (MW) & time (hr) &  rate (MW/h) \\
    \midrule
    Natural gas & 20 & 33.4 & 1-5 & 3 & 2.5 \\
    Hydrogen & 1 & 30 & $10^{-5}$-1 & N/A & N/A\\
    \bottomrule
    \toprule
    Storage & \# of & Storage cost & Capacity & Min-max charging/ & Min charging/ \\
    type & units & (\$/ton) & & discharging & discharging\\
    & & & & power (MW) & time (hr)\\
    \midrule
    ESS & 1 & N/A & 20 (MWh) & 0.8-4 & 5 \\
    HSC & 1 & 10,000 & 10 (ton) & N/A & N/A \\
    \bottomrule
    \end{tabular}
    \caption{Characteristics of fossil based dispatchable units \cite{khodaei2013microgrid}, ESS, and hydrogen storage.} \label{table_parameters}
\end{table}

\begin{figure}[ht!]
    \centering
    \includegraphics[width=0.8\textwidth]{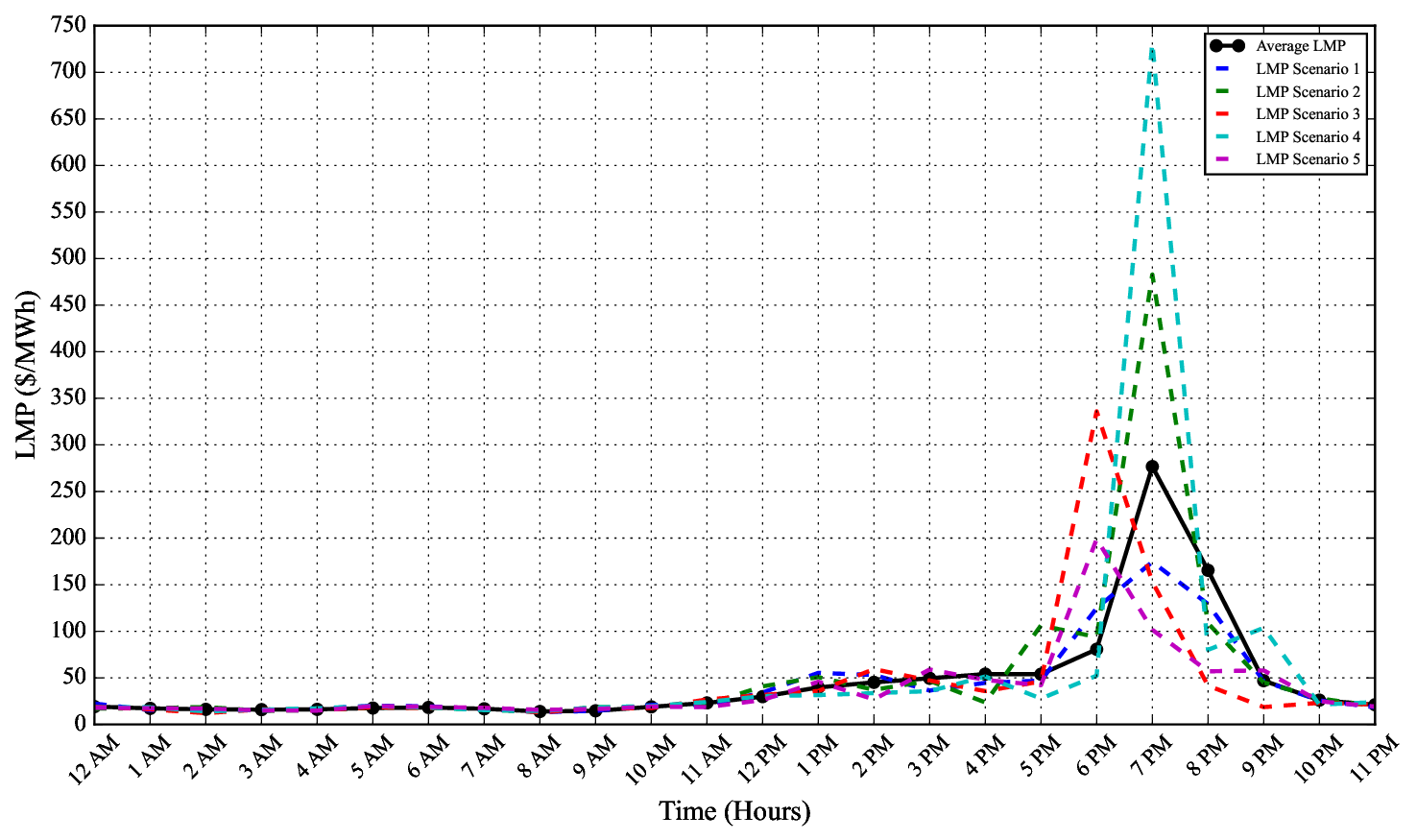}
    \vspace{-1em}
    \caption{The output electricity market prices of five scenarios considered.}\label{fig_lmp}
\end{figure}

\begin{figure}[ht!]
    \centering
    \includegraphics[width=0.8\textwidth]{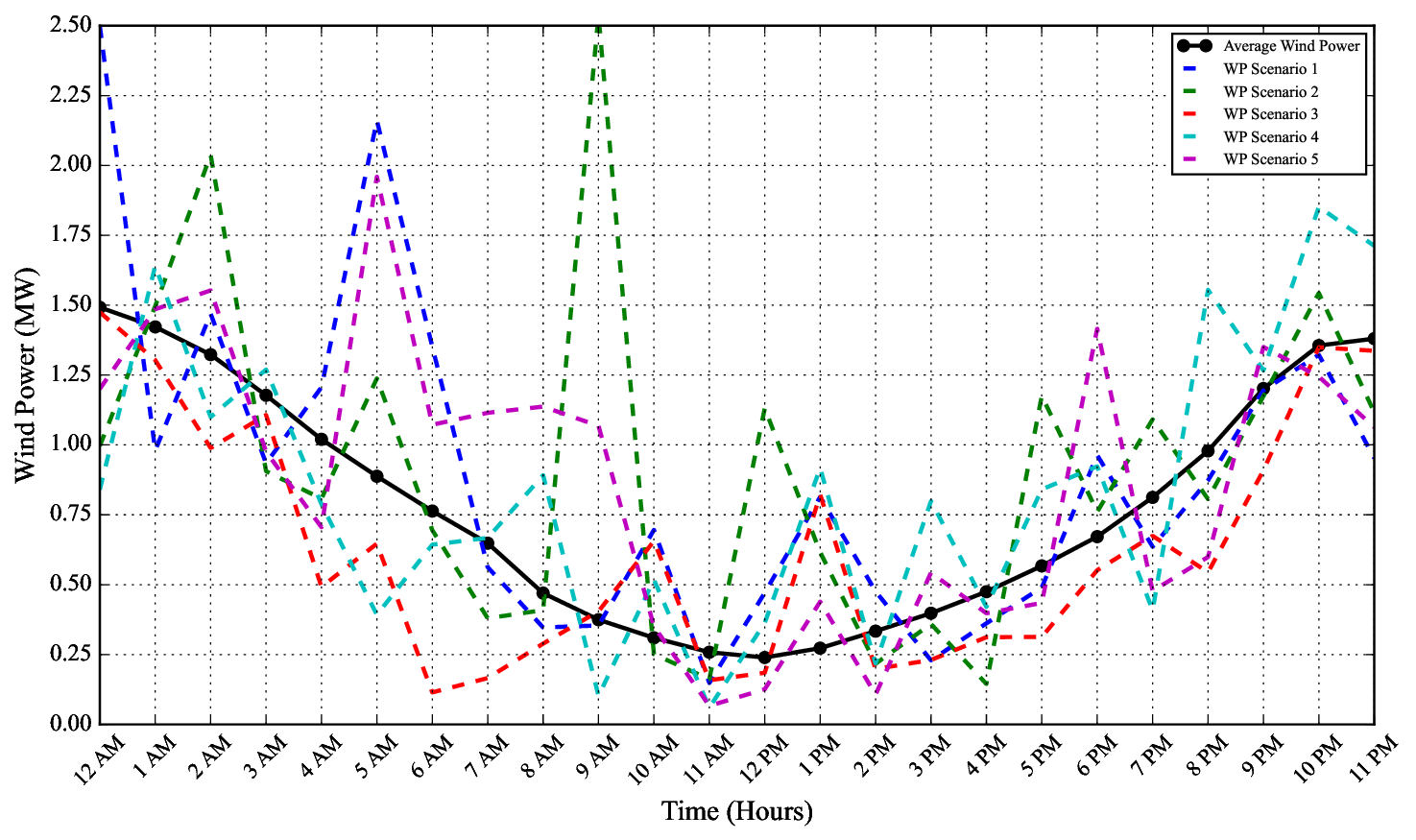}
    \vspace{-1em}
    \caption{The output power of wind turbines of five scenarios considered.}\label{fig_wind}
\end{figure}

\begin{figure}[ht!]
    \centering
    \includegraphics[width=0.8\textwidth]{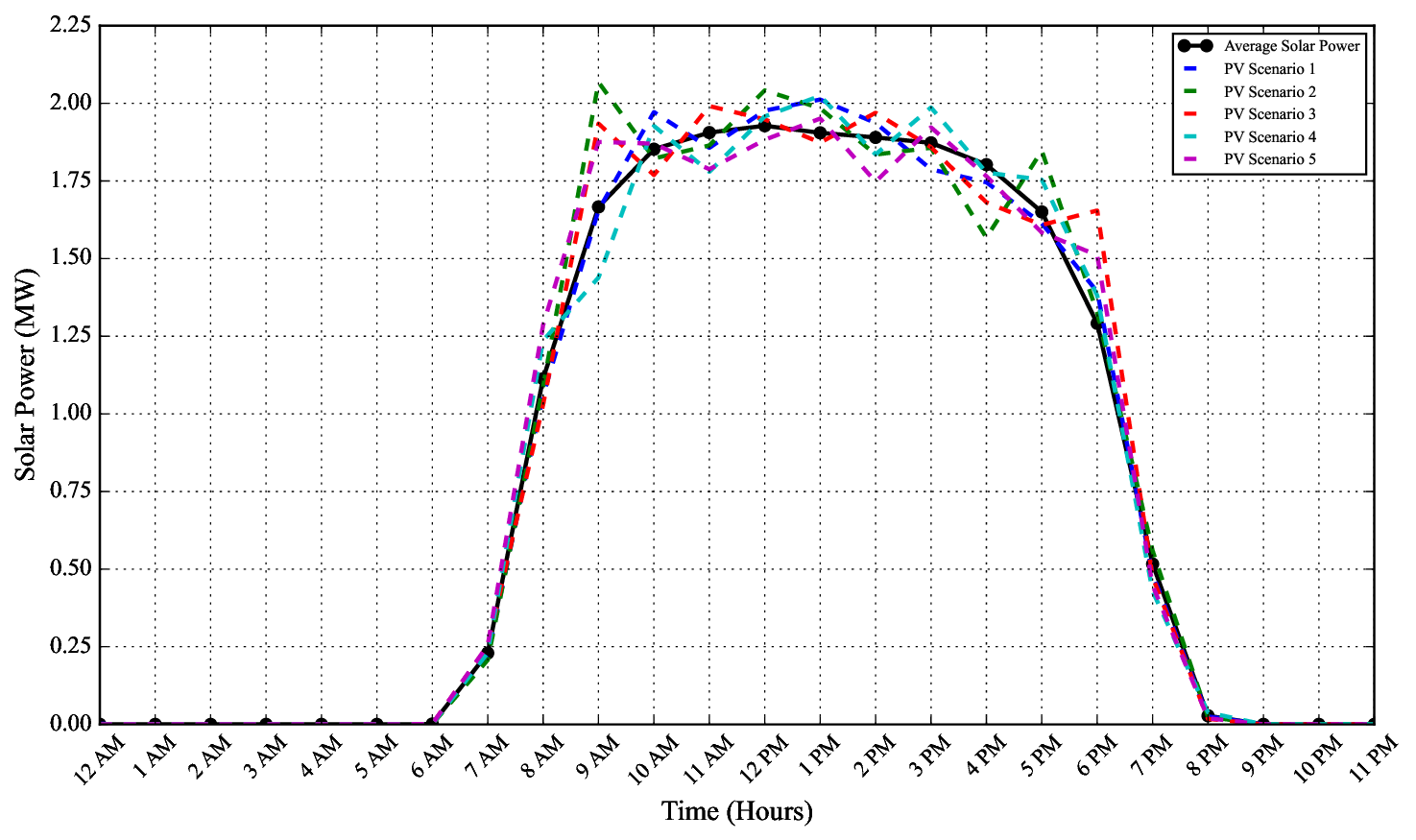}
    \vspace{-1em}
    \caption{The output power of photovoltaic panels of five scenarios considered.}\label{fig_solar}
\end{figure}

The electricity market price as well as solar and wind energy generation are intermittent in nature. To account for these uncertainties, we synthesize five scenarios for each uncertain parameter using Monte Carlo simulations on a log-normal distribution to prevent getting any negative parameter values. These scenarios are generated based on publicly available Texas grid data from ERCOT in August 2024. For locational marginal pricing (LMP), we randomly selected Bus TC-KO without loss of generality. Table \ref{table_parameters_scenario} summarizes the probability of scenarios for each uncertain parameter. Note that there still will be $5^3 = 125$ combinations of these scenarios. To mitigate this, we further perform a scenario reduction to obtain 5 representative combinations using a probability distance algorithm based on the Kantorovich distance \cite{gholami2016microgrid} (see Table \ref{table_scenario}). Figures \ref{fig_lmp} through \ref{fig_solar} illustrate the mean and average electricity market price, wind power, and solar power for the 5 combined scenarios.
\begin{table}[ht]
    \centering 
    \begin{tabular}{cccc}
    \toprule
    Scenario number used & & Probability of uncertain parameter & \\
    in Table \ref{table_scenario} & LMP & WP & PV \\
    \midrule
    1 & 64.50\% & 13.48\% & 4.86\%\\
    2 & 2.92\% & 1.52\% & 30.31\% \\
    3 & 3.26\% & 64.91\% & 15.81\% \\
    4 & 1.29\% & 6.18\%& 32.96\% \\
    5 & 28.03\% & 13.91\%& 16.06\% \\
    \bottomrule
    \end{tabular}
    \caption{Probability of generated scenarios for LMP, wind power, and PV.} \label{table_parameters_scenario}
\end{table}

\begin{table}[ht]
    \centering 
    \begin{tabular}{ccccc}
    \toprule
    Combined & Selected LMP  &  Selected WP  &  Selected PV  & Probability\\
    scenario & scenario & scenario & scenario & \\
    \midrule
    1 & 5 & 3 & 2 & 18.21\%\\
    2 & 1 & 3 & 4 & 45.56\%\\
    3 & 1 & 1 & 3 & 4.54\%\\
    4 & 5 & 3 & 3 & 9.50\%\\
    5 & 1 & 3 & 5 & 22.19\%\\
    \bottomrule
    \end{tabular}
    \caption{Probability of generated scenarios (after scenario reduction) for the microgrid problem after scenario reduction.} \label{table_scenario}
\end{table}

In this chapter, we will study the optimal scheduling problem under both the grid-connected and islanded modes. We will also closely examine two different scenarios in the grid-connected mode to draw insights regarding how the microgrid should be scheduled and operated. 

\subsection{Grid-connected Mode}

\begin{figure}[ht]
    \centering
    \includegraphics[width=\textwidth]{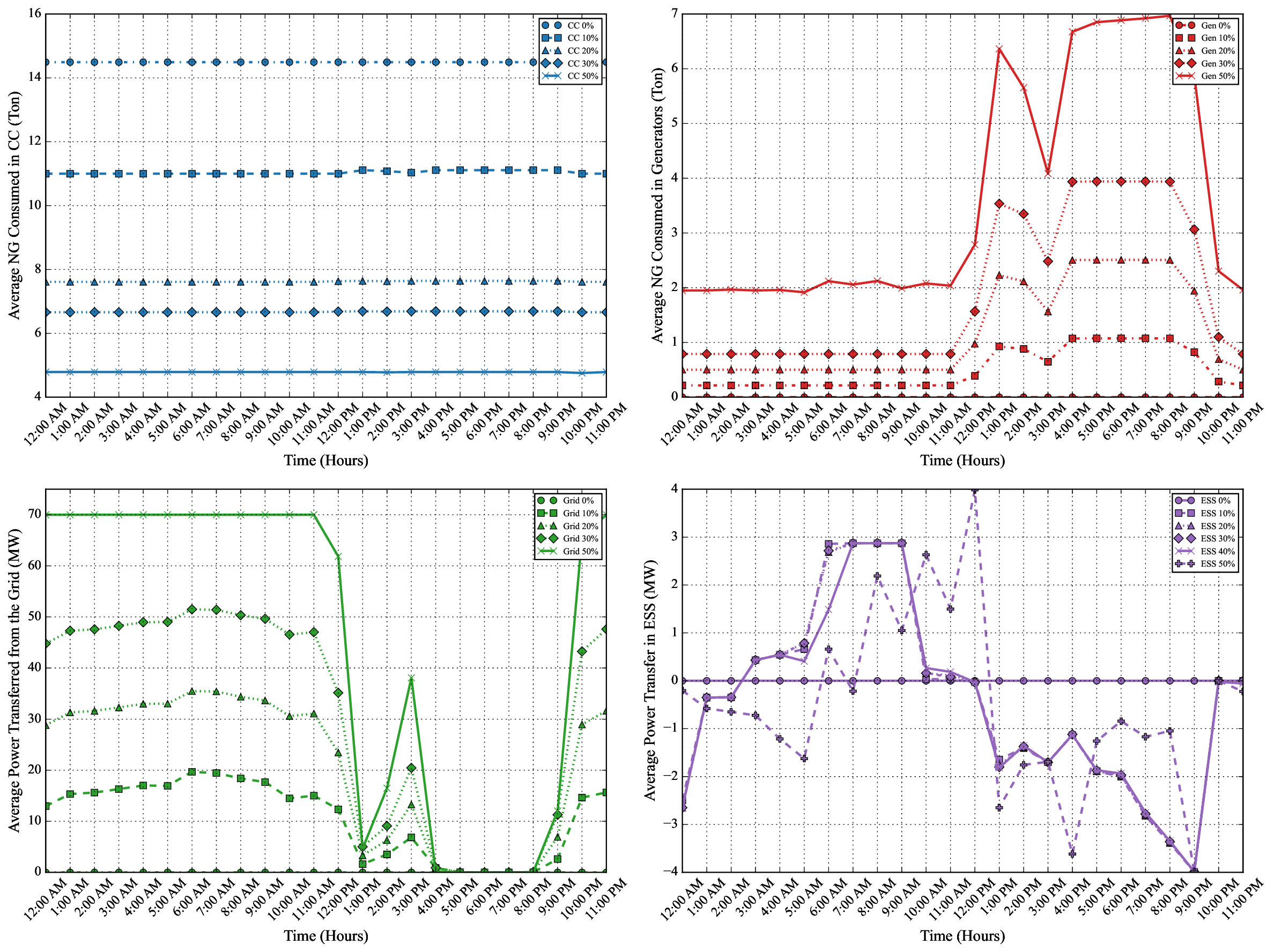}
    \vspace{-2em}
    \caption{Average natural gas (NG) consumption for conventional crackers, natural gas consumption for dispatchable generators, power usage from the main grid, and ESS charging/discharging status for grid-connected mode under different degrees of electrification.} \label{fig_case1_NG}
\end{figure}

\begin{figure}[ht]
    \centering
    \includegraphics[width=\textwidth]{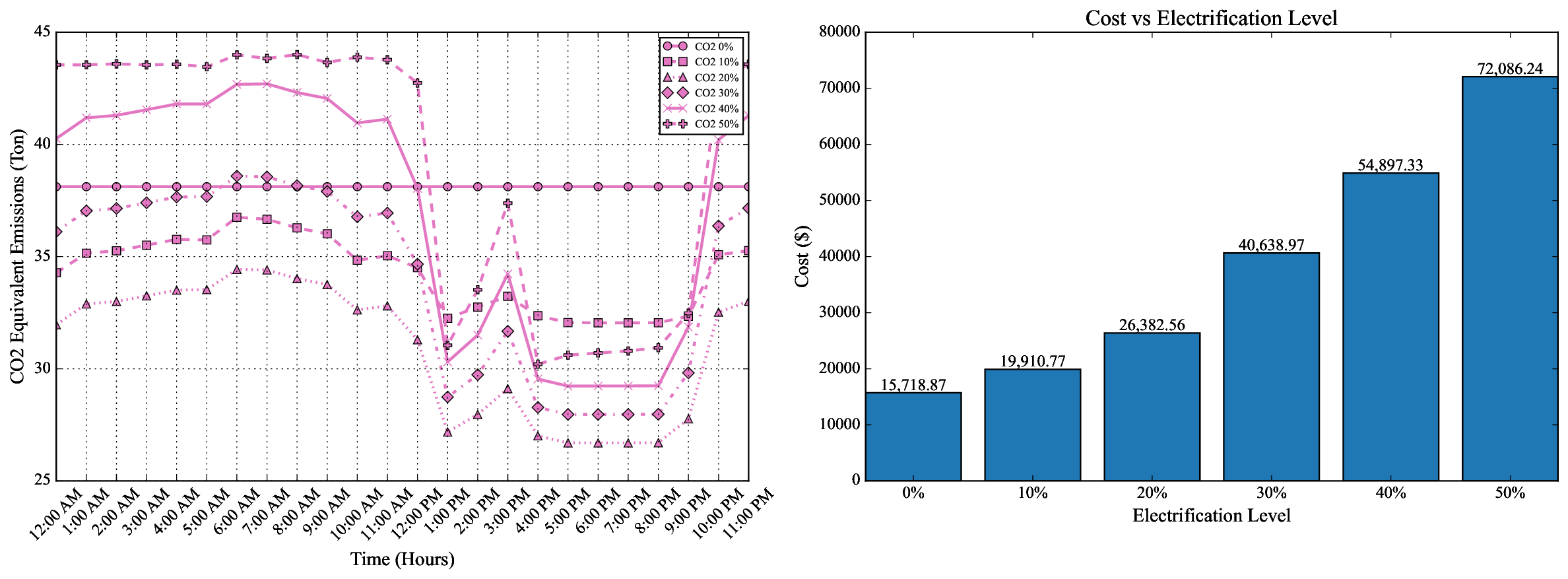}
    \vspace{-2em}
    \caption{Grid-connected mode: average \ce{CO2} equivalent emissions and expected total cost at different electrification levels.} \label{fig_case1_co2}
\end{figure}

In this case, we consider the (partially) electrified ethylene plant to be operated in the grid-connected mode for a 24-hour horizon. The ESS is considered to be half-charged at $t=0$. Specifically, since we envision a gradual electrification progress in chemical plants, we minimize the total cost of Equation \eqref{eqn_obj} while considering 0\%, 10\%, 20\%, 30\%, 40\%, and 50\% of the ethylene is produced in electrified crakers. In Figure \ref{fig_case1_NG}, we show the amount of natural gas consumed by conventional crackers and fossil-based generators, as well as power usage from the grid and ESS charging status. In Figure \ref{fig_case1_co2}, we show the hourly \ce{CO2} equivalent emissions and the objective function values corresponding to different electrification levels.

The results show that, as the plant gets more electrified, the dependency of the microgrid on the main grid increases, and the expected total cost increases as well. At 50\% electrification, the expected total cost is 4.6 times that with 0\% electrification. On the other hand, the average hourly \ce{CO2} equivalent emissions corresponding to 0\% through 50\% electrification are 38.116, 34.303, 30.951, 34.093, 37.234, and 39.219 tons, respectively. This is an interesting observation because increasing the electrification level does not always lead to a reduction in carbon intensity. Rather, the optimal degree of electrification lies somewhere 10\% to 30\% from a sustainability point of view. For example, at 20\% electrification, the average hourly \ce{CO2} equivalent emissions at the optimal solution (from a total expected cost perspective) is 18.8\% lower than the conventional case with 0\% electrification. Meanwhile, the optimal expected total cost is only 68\% higher than the conventional case and is significantly lower than the 50\% electrification case.

By cross-referencing Figures \ref{fig_case1_co2} and \ref{fig_case1_NG}, we believe the reason behind this is the following. Between 0\% to 50\% electrification, the average natural gas consumption by conventional crackers declines as expected. Nevertheless, the rate of decrease is the greatest between 0\% to 20\% electrification and then diminishes beyond 20\% electrification. On the other hand, the natural gas consumed by dispatchable generators shows a steady increase as the level of electrification increases. In fact, the average hourly natural gas consumption by conventional crackers as well as by fossil-based generators for 0\%, 10\%, 20\%, 30\%, 40\%, and 50\% electrification are 14.493, 11.541, 8.814, 8.549, 8.284, and 8.429 ton, respectively. As a result, the overall natural gas consumption flattens out at and after 20\% of electrification. Meanwhile, when renewable electricity generation technologies (e.g., fuel cells) are not yet mature enough to be cost-competitive compared to conventional fossil-based generators, as the power demand from electrified cracker units increases, the microgrid must produce more electricity using its local generators, especially during peak hours when the LMP is high. However, due to efficiency losses, this will incur more energy demand than directly using natural gas as fuel for conventional crackers, thereby incurring more costs and higher carbon intensity.

Another observation we can draw is that, when the electricity price is high, the power withdrawn from the main grid decreases drastically to 0. In turn, the local generators will pick up the power demand, and the ESS will start discharging to compensate for the missing demand. During non-peak hours, the ESS will take advantage of the low electricity price to charge itself and prepare for the anticipated demand during peak hours. When the electricity price is high, the \ce{CO2} equivalent emissions decrease significantly, which is consistent with the observation that local renewable electricity generation technologies are participating more heavily in operating the crackers. Meanwhile, we remark that in order to achieve holistic decarbonization via electrification, it is important for both power systems and chemical plant stakeholders to decarbonize at a similar pace. For example, if clean energy sources contribute a higher proportion to ERCOT's energy profile, then a higher degree of electrification in the ethylene plant will translate to lower carbon emissions (which is not what we currently observe in this case study).

\begin{figure}[ht]
    \centering
    \includegraphics[width=\textwidth]{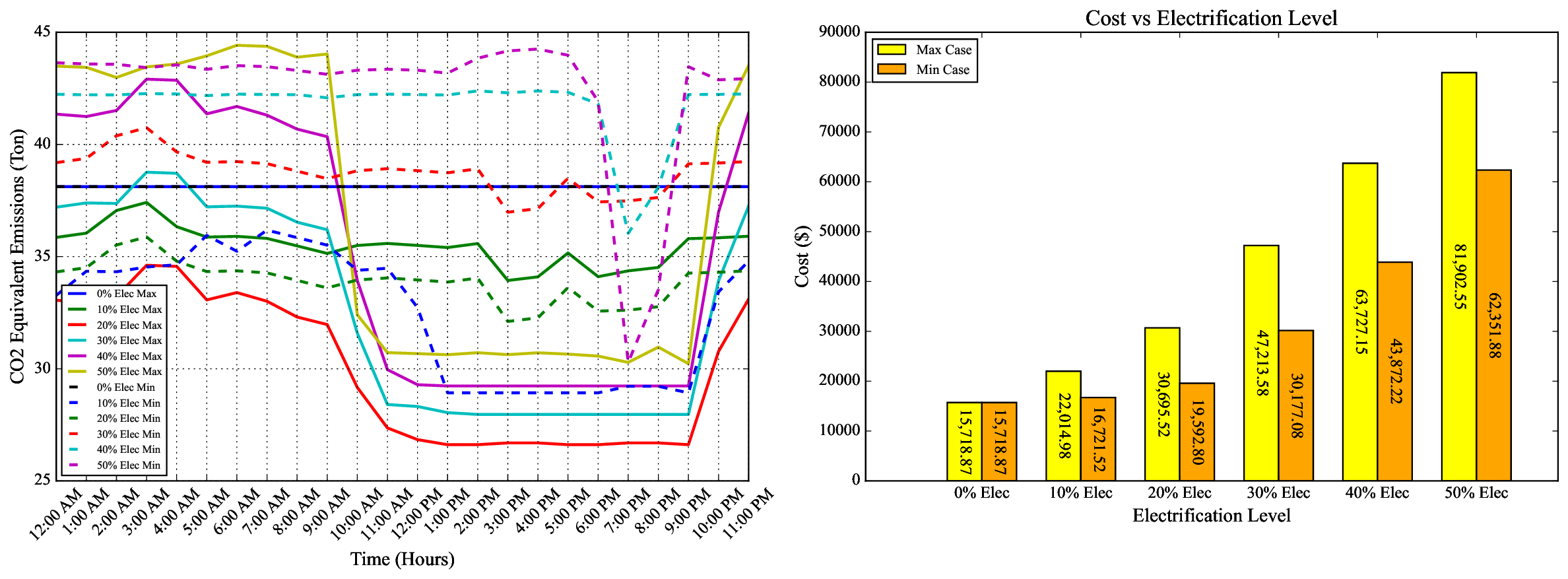}
    \vspace{-2em}
    \caption{Comparing \ce{CO2} equivalent emissions and expected total cost under different electricity prices.} \label{fig_comparison}
\end{figure} 

Next, we investigate the impact of electricity prices from the grid on electrification progress. As shown in Figure \ref{fig_comparison}, we compare the \ce{CO2} equivalent emissions and total expected cost of the optimal solution when using the lowest hourly electricity market price among the five scenarios as well as when using the highest. As we can see, with higher electricity prices, the overall carbon intensity is reduced as the plant prefers to generate its power locally. This is because the microgrid, which harnesses the full decarbonization potential (see Figure \ref{fig_superstructure}) by considering various renewable options, is generally more renewable compared to the main grid. On the other hand, when using the lowest hourly electricity price, the plant prefers to be more grid-dependent. And the average hourly \ce{CO2} equivalent emissions are 38.116, 32.524, 33.928, 38.798, 41.792, and 42.454 tons for electrification levels of 0\% to 50\%, respectively, which are higher than the carbon intensity values shown in Figure \ref{fig_case1_co2}. Clearly, there is a trade-off between process economics (total cost) and sustainability (carbon intensity) attributed to the market price and carbon intensity of power from the main grid. From this illustrative example, we conclude that the decarbonization of chemical manufacturing via electrification must be accompanied by the simultaneous decarbonization of power systems.

\subsection{Islanded Mode}

\begin{figure}[ht]
    \centering
    \includegraphics[width=\textwidth]{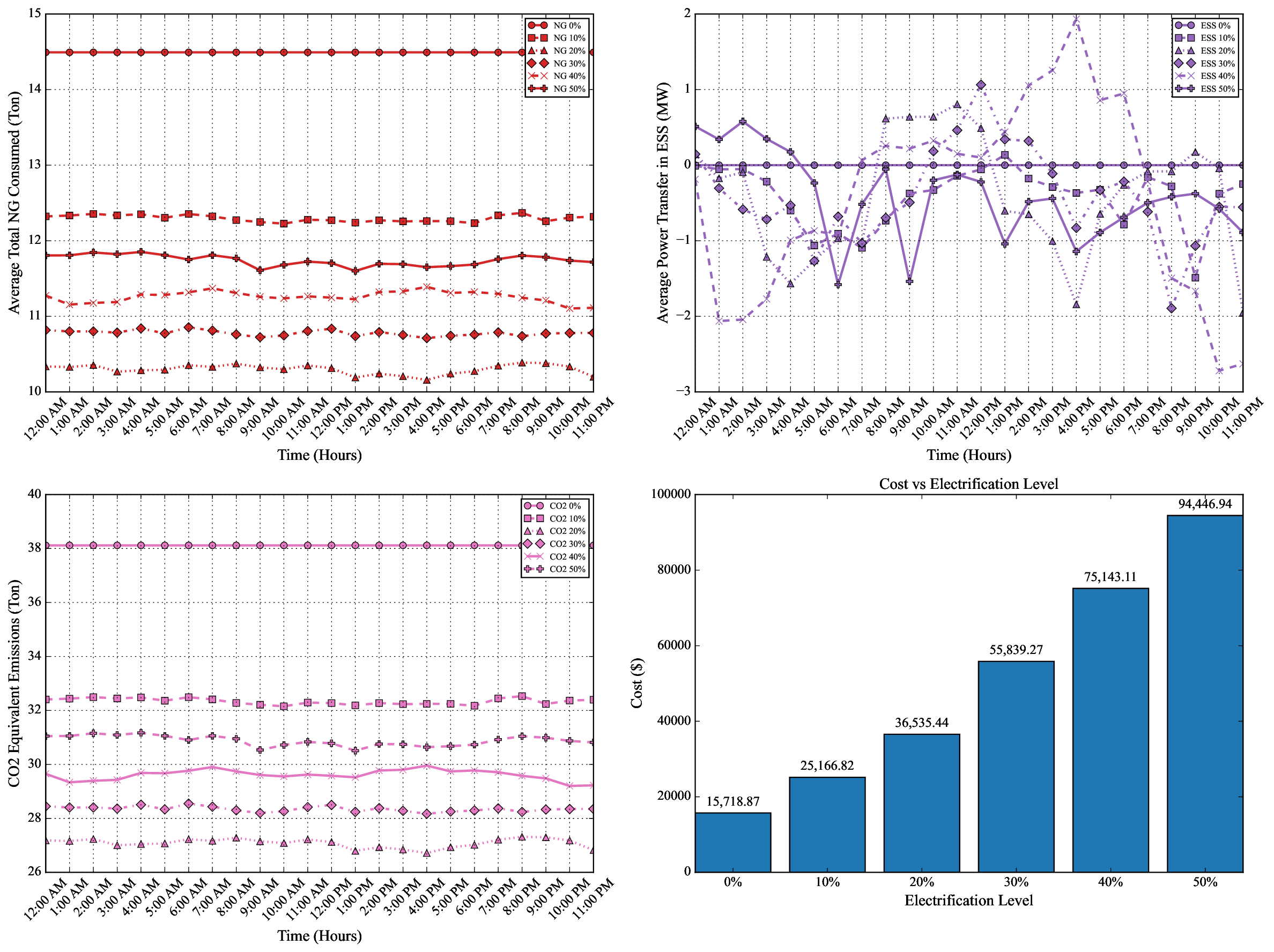}
    \vspace{-2em}
    \caption{Islanded mode: Average hourly total natural gas consumption, energy storage system charging/discharging status, average \ce{CO2} equivalent emissions, and expected total cost.} \label{fig_case2}
\end{figure}

Lastly, we consider the case where the ethylene plant is responsible for satisfying all of its energy demand locally, which could occur during a power system outage. We assume that the energy storage system is initially half-charged. In Figure \ref{fig_case2}, we show the total hourly natural gas consumption (accounting for natural gas used for conventional crackers and dispatchable generators), energy storage charging/discharging status, \ce{CO2} equivalent emissions, and total expected cost for 0\%, 10\%, 20\%, 30\%, 40\%, and 50\% electrification. In terms of \ce{CO2} equivalent emissions, we find that, even though the carbon intensity is much lower when the plant is operated in islanded mode, the optimal degree of electrification under the current technology status is still around 20\% (i.e., roughly 1 electrified cracker and 4 conventional crackers). By increasing ESS and \ce{H2} storage capacities, we may debottleneck the decarbonization barrier and achieve lower emissions with a higher degree of electrification.

\section{Conclusion} \label{sec_conclusion}

Steam cracking is one of the most important and energy-intensive chemical processes that exhibit great opportunities for decarbonization. In this book chapter, we develop a microgrid superstructure for the steam cracking process considering diverse energy sources and storage systems. To obtain the minimum energy requirement for conventional and electrified steam cracking, we formulate a DAE-based constrained optimization problem and solve it as an NLP by discretization using the orthogonal collocation method in \texttt{pyomo.dae}. Furthermore, we present a single-stage scenario-based MILP formulation for optimal microgrid scheduling. Using actual VRE profiles and electricity market prices, we study an illustrative case study of microgrid scheduling under grid-connected and islanded modes. The carbon intensity and expected total cost under different decarbonization levels are determined along with the corresponding electricity and natural gas usage. The results show that, given the current status of the power grid and renewable energy generation technologies, the process economics and sustainability of electrified steam cracking do not always favor higher decarbonization levels. To achieve true decarbonization of olefins production, electricity from the main grid must be cleaner and cheaper, and energy storage (ESS and \ce{H2} storage) costs per unit stored must go down. Furthermore, it is important for both chemical and power systems stakeholders must seamlessly coordinate with each other to pursue joint optimization in operation. 

\section*{Acknowledgements}
This work is supported by the U.S. National Science Foundation under Award 2343072 and the start-up fund of the College of Engineering, Architecture, and Technology at Oklahoma State University.

\end{document}